\documentstyle[12pt,amssymb]
{article}

\title{Quantum  ``az+b'' group
at roots of unity: unitary representations }

\author{Ma{\l }gorzata Rowicka\thanks{Supported by KBN grant 
No 5 PO3A 045 20}\\ 
Institute of Mathematics, Polish Academy of Sciences,\\
\'Sniadeckich 8, 00-950 Warszawa, Poland\\
email: rowicka@fuw.edu.pl}

\begin{document}

\maketitle
\newcommand{\qed}{$\Box$}
\newcommand{\mqed}{\nopagebreak\centerline{\hfill
\raisebox{3.5ex}[0ex][0ex]{$\Box$}}\\}
\setlength{\parindent}{0in}
\setlength{\parskip}{4mm}
\newcommand{\faz}{{\rm Phase\  }}
\newcommand{\ci}{continuous}
\newcommand{\ru}{ unitary representation}
\newcommand{\Qa}{ Quantum 'az+b' group }
\newcommand{\qa}{ quantum 'az+b' group }
\newcommand{\MU}{ multiplicative unitary}
\newcommand{\slw}{S.L. Woronowicz}
\newcommand{\tr}{\Delta}
\newcommand{\cH}{{\cal H}}
\newcommand{\K}{{\cal K}}
\newcommand{\cK}{{\cal K}}
\newcommand{\sH}{_{\cal H}}
\newcommand{\sK}{_{\cal K}}
\newcommand{\h}[1]{\hat{#1}}
\newcommand{\rf}[1]{{\rm (\ref{#1})}}
\newcommand{\ch}{ {\cal C}(H)}
\newcommand{\chn}{ {\cal C}(H)^N}
\newcommand{\ov}{\overline}
\newcommand{\fh}{F_{\hbar}} 
\newcommand{\hb}{\hbar}
\newcommand{\vt}{V_{\theta}} 
\newcommand{\za}{-\!\!\circ} 
\newcommand{\R}{{\mathbb R}}
\newcommand{\C}{{\mathbb C}}
\newcommand{\Z}{{\mathbb Z}}
\newcommand{\ro}{\rho}
\newcommand{\si}{\sigma}
\newcommand{\be}{\beta}
\newcommand{\de}{\delta}
\newcommand{\ga}{\gamma}
\newcommand{\ta}{\tau}
\newcommand{\N}{{\mathbb N}}
\newcommand{\M}{{\rm M}}
\newcommand{\B}{{\rm B}}
\newcommand{\Lin}{{\rm L}}
\newcommand{\pod}{{\rm d}}
\newcommand{\iz}{\cong}
\newcommand{\eps}{\epsilon}
\newcommand{\fil}{\varphi}
\newcommand{\la}{\lambda}
\newcommand{\Ci}{C_{\infty}}
\newcommand{\Cir}{C_{\infty}(\R)}
\newcommand{\eh}{e^{\frac{i\hbar}{2}} }
\newcommand{\ehm}{e^{-\frac{i\hbar}{2}} }
\newcommand{\Cgr}{C^{\infty}(\R)}
\newcommand{\Cg}{C^{\infty}}
\newcommand{\Cor}{C_{\rm o}(\R)}
\newcommand{\Co}{C_{\rm o}}
\newcommand{\Cog}{C_{\rm bounded}}
\newcommand{\Cogr}{C_{\rm b}(\R)}
\newcommand{\Lk}{L^{2}}
\newcommand{\Lkr}{L^{2}(\R)}
\newcommand{\po}{\hat{p}}
\newcommand{\qo}{\hat{q}}
\newcommand{\xo}{\hat{x}}
\newcommand{\csta}{$C^{*}$-}
\newcommand{\cstal}{$C^{*}$-algebra}
\newcommand{\te}{\otimes}
\newcommand{\ad}{{\rm ad}}
\newcommand{\id}{{\rm id}}
\newcommand{\Mor}{{\rm Mor}}
\newcommand{\Rep}{{\rm Rep}}
\newcommand{\spe}{{\rm Sp }}
\newcommand{\sign}{{\rm sign }\;}
\newcommand{\whe}{\hspace*{5mm}\mbox{\rm where}\hspace{5mm}}
\newcommand{\mand}{\hspace*{5mm} {\rm and} \hspace{5mm}}
\newcommand{\moraz}{\hspace*{5mm} {\rm and} \hspace{5mm}}
\newcommand{\af}{\hspace*{1mm} {\bf \eta} \hspace{1mm}}
\newcommand{\od}{\hspace*{5mm}}
\newcommand{\fu}{{\cal F}}
\newcommand{\fuod}{{\cal F}^{-1}}
\newcommand{\mlot}{\mbox{$\hspace{.5mm}\bigcirc\hspace{-3.7mm}
\raisebox{-.7mm}{$\top$}\hspace{1mm}$}}
\newcommand{\Zak}{\mbox{$-\hspace{-2pt}\comp\,$}}
\newcommand{\dow}{{\bf Proof: }}
\newcommand{\ut}{\cong}
\newcommand{\Sp}{{\rm sp}}
\newcommand{\infi}{\infty}
\newcommand{\tend}{\rightarrow}
\newcommand{\impl}{\Rightarrow}
\newcommand{\limn}{\lim_{n\rightarrow\infty}}
\newcommand{\limk}{\lim_{k\rightarrow\infty}}
\newcommand{\limt}{\lim_{t\rightarrow\infty}}
\newcommand{\limx}{\lim_{|x|\rightarrow\infty}}
\newcommand{\lime}{\lim_{\eps\rightarrow 0}}
\newcommand{\Lj}{L^{1}}
\newcommand{\Lp}{L^{P}}
\newcommand{\Ls}{L^{S}}
\newcommand{\Li}{L^{\infty}}
\newcommand{\ljr}{{\cal L}^{1}(\R)}
\newcommand{\lk}{ l^{2}}
\newcommand{\lpr}{{\cal L}^{P}(\R)}
\newcommand{\lir}{{\cal L}^{\infty}(\R)}
\newcommand{\Ljr}{L^{1}(\R)}
\newcommand{\Lpr}{L^{P}(\R)}
\newcommand{\Lsr}{L^{S}(\R)}
\newcommand{\Lir}{L^{\infty}(\R)}
\newcommand{\lkn}{ l^2(\N)}
\newcommand{\lp}{{\cal L}^{P}}
\newcommand{\li}{{\cal L}^{\infty}}
\newcommand{\Ljx}{L^{1}(\X)}
\newcommand{\Lkx}{L^{2}(\X)}
\newcommand{\Lpx}{L^{p}(\X)}
\newcommand{\Lqx}{L^{q}(\X)}
\newcommand{\Lix}{L^{\infty}(\X)}
\newcommand{\ljx}{{\cal L}^{1}(\X)}
\newcommand{\lkx}{{\cal L}^{2}(\X)}
\newcommand{\lpx}{{\cal L}^{p}(\X)}
\newcommand{\lix}{{\cal L}^{\infty}(\X)}
\newcommand{\ilkC}{C_{\infty}(\R)\otimes_{C}C_{\infty}(\R)}
\newcommand{\ilk}{\otimes_{\cal C}}
\newcommand{\zw}{CB(\Lkr)}
\newcommand{\zwg}{CB\left(\,L^{2}(G)\,\right)}
\newcommand{\ogr}{B(L^2(\R))}
\newcommand{\hs}{HS(H)}
\newcommand{\spl}{\star}
\newcommand{\milk}{$(i_{1},i_{2},\zw)$}
\newcommand{\mczw}{{\rm Mor(\Cir,\zw)}}
\newcommand{\pra}{\mbox{${\rm Proj_1}$}}
\newcommand{\prb}{\mbox{${\rm Proj_2}$}}
\newcommand{\Morc}{{\rm Mor}_{\katC}}
\newcommand{\recogr}{{\rm Rep(\Cir,\ogr)}}
\newcommand{\mi}{\hspace*{3mm} {\rm and} \hspace{3mm}}
\newcommand{\mor}{\hspace*{5mm} {\rm or} \hspace{5mm}}
\newcommand{\dla}{\hspace*{5mm} {\rm for} \hspace{5mm}}
\newcommand{\ja}{j_{1}}
\newcommand{\ib}{i_{2}}
\newcommand{\ia}{i_{1}}
\newcommand{\jb}{j_{2}}
\newcommand{\jc}{j_{3}}
\newcommand{\sz}{{\cal S}(\R)}
\newcommand{\gwhom}{$^*$-homomorfizm}
\newcommand{\gw}{$^*$}
\newcommand{\heps}{h_{\eps}} 
\newcommand{\czn}{\Co}
\newcommand{\bh}{B(H) }
\newcommand{\dlad}{for any }
\newcommand{\dlak}{for each }
\newcommand{\Gn}{{\cal G}_n}
\newcommand{\G }{{\rm G}_n}
\newcommand{\rt}{{\bf T}^1}
\newcommand{\jd}{{\bf T}^2}
\newcommand{\Tau}{{\cal T}}


\newcommand{\bfa}{\begin{fakt}}\newcommand{\efa}{\end{fakt}}
\newcommand{\ble}{\begin{lem}}\newcommand{\ele}{\end{lem}}
\newcommand{\bst}{\begin{stw}}\newcommand{\est}{\end{stw}}
\newcommand{\bde}{\begin{defi}}\newcommand{\ede}{\end{defi}}
\newcommand{\bwn}{\begin{wn}}\newcommand{\ewn}{\end{wn}}
\newcommand{\buw}{\begin{uwaga}}\newcommand{\euw}{\end{uwaga}}
\newcommand{\bdy}{\begin{dygresja}}\newcommand{\edy}{\end{dygresja}}
\newcommand{\bwa}{\begin{warning}}\newcommand{\ewa}{\end{warning}}
\newcommand{\bpr}{\begin{przy}}\newcommand{\epr}{\end{przy}}
\newcommand{\btw}{\begin{tw}}\newcommand{\etw}{\end{tw}}
\newcommand{\beq}{\begin{equation}}\newcommand{\eeq}{\end{equation}}
\newcommand{\bit}{\begin{itemize}}\newcommand{\eit}{\end{itemize}}
\newcommand{\bq}{\begin{quote}}\newcommand{\eq}{\end{quote}}
\newcommand{\ba}{\begin{array}}\newcommand{\ea}{\end{array}}


\newtheorem{defi}{Definition}[section]
\newtheorem{wn}[defi]{Observation}
\newtheorem{tw}[defi]{Theorem}
\newtheorem{lem}[defi]{Lemma}
\newtheorem{fakt}[defi]{Corollary}
\newtheorem{stw}[defi]{Proposition}
\newtheorem{przy}[defi]{Example}
\newtheorem{uwaga}[defi]{Remark}

\begin{abstract} 
All unitary representations of the  
 quantum ``az+b'' group are found . It turns out that this quantum group is 
self dual i.e. all unitary representations are  
 'numbered' by elements of the same group. Moreover, 
 the  formula for all unitary representations 
involving the quantum exponential function is proven.
\end{abstract}
{\bf keywords: } \cstal -- crossed product, quantum group\\
{\bf MSC-class: }{ 20G42 (Primary), 47L65, 81R15 (Secondary).}

\section{Introduction}
Locally compact quantum groups are nowadays studied extensively 
by many scientists \cite{ks}.Although at the moment there 
is no commonly accepted definition of a locally compact quantum group,
there are promising approaches and interesting examples have been worked out. 
One of the most remarkable  ones is the quantum ``az+b''
 group constructed by S.L. Woronowicz in \cite{az+b}.
According to the recent computation by A. van Daele \cite{nowywandal}, 
this group is  an example of an interesting
phenomenon  foreseen by 
Vaes and Kustermans in \cite{vaesk}.
In this paper we study the \qa      
  from the point of view of unitary representations and duality theory.
The aim of this paper is to derive a formula for all unitary 
representations of the \qa.

In Section  \ref{qa} 
 we recall relevant information on the quantum  ``az+b'' group. 
We discuss in details the structure of the  \cstal $\;$  
 of all  continuous functions vanishing at infinity on 
 the quantum  ``az+b'' group. 
We introduce also the  corresponding  $C^*$- and  $W^*$-dynamical 
systems. 
 We prove a couple of corollaries we use later on 
 to find all unitary representations of the quantum  ``az+b'' group. 
The formula for this representations will be found and  
proved in  Section \ref{ost}. 
The important theorem we use 
 was proved in our previous paper \cite{paper3}, where we investigated  
 unitary representations 
 of some braided quantum group related to the quantum ``az+b''. 

 We  use methods  similar to those introduced by 
S.L. Woronowicz in the case of the quantum  $E(2)$ group  \cite{oe2} and 
then applied by us to  
 giving the formula for all unitary  representations  
 of the quantum $''ax+b''$ group in  \cite{paper2} and 
in Chapter 3 of   \cite{phd}.

In the next Section we will fix the  notation.

\section{Notation}
\label{not}

We consider only concrete \cstal s, i.e. embedded into  \cstal $\;$
\ of all bounded operators acting  on Hilbert space $\cH$, denoted by $B(\cH)$.
The \cstal $\;$
\  of all compact  operators acting  on  $\cH$ will be  denoted by $CB(\cH)$.
All algebras we consider are separable with the exception of multiplier algebras (see definition of multiplier algebra below).

Let $A$ be \cstal. Then $M(A)$ will denote the
{\em multiplier algebra} 
of $A$, i.e.
\[M(A)=\{m\in B(\cH):\od ma,am\in A\  \ \mbox{\dlad} \ a\in A\}\ .\]
Observe that $A$ is an  ideal in $M(A)$. If   $A$ is a unital 
 \cstal\ , then $A=M(A)$, in general case $A\subset M(A)$.
For example the multiplier algebra of $CB(\cH)$ is the algebra  $B(\cH)$ and 
the multiplier algebra of \csta algebra $\Cir$ 
of all continuous vanishing at infinity functions on $\R$ is the algebra of all continuous bounded functions on $\R$ denoted by $C_{\scriptstyle bounded}(\R)$.
The natural topology on $M(A)$ is the strict topology, 
i.e. we say that a sequence $(m_n)_{n\in\N}$ of $m_n\in M(A)$ 
converges strictly to 0 if for every $a\in A$, we have 
$||m_n a||\rightarrow 0$ and $||a m_n ||\rightarrow 0$, when 
$n\rightarrow +\infty$. 
Whenever we will consider continuous maps from or into 
$M(A)$, we will mean this topology. 

For any \cstal s $A$ and $B$, we will say that $\phi$ is a morphism 
and write $\phi\in \Mor (A,B)$ if  $\phi$ is a * - algebra
homomorphism acting from $A$ into $M(B)$ and such that 
$\phi (A)B$ is dense  in $B$.
Any $\phi\in \Mor (A,B)$ admits unique extension to a * - algebra
homomorphism acting from $M(A)$ into $M(B)$. For any $S\in M(A)$,
operator $\phi(S)$ is given by 
\[\phi(S)(\phi(a)b)=\phi(Ta)b\ ,\]
where $a\in A$ and $b\in B$.

For any closed  operator $T$ acting on   $\cH$ 
 we define its  {\em $z$-transform} by
\[z_T=T(I+T^*T)^{-\frac{1}{2}}\ .\]
\label{ztransf}
Observe that   $z_T\in \B(\cH)$ and  $||z_T||\leq 1$.
Moreover, one can recover  $T$ from $z_T$
\[T=z_T(I-z_T^*z_T)^{-\frac{1}{2}}\ .\]

A closed operator $T$ acting on $A$ is {\em affiliated} with a \cstal $\;$ \ $A$ 
if and only if $z_T\in \M (A)$ and  $(I-z_T^*z_T)A$ 
is dense in $A$. 
 A set of all elements affiliated with $A$ is denoted by  
$A^{\af}$\index{$A^{\af}\;\;$\}.
If  $A$ is a unital \cstal , then  $A^{\af}=M(A)=A$, 
 in general case 
\[A\subset M(A)=\{T\af A\;:\;||T||<\infty \}\subset A^{\af} \ .\]
The set of all elements affiliated with  $\Cir$ is the set of all  
continuous functions on 
real line $C(\R)$, and a 
 set of all  elements affiliated  with \cstal $\;$\ $CB (\cH)$ is a set 
 of all closed operators  $\ch$. 
This last example shows that a product and a sum of two
 elements affiliated with  $A$ may not  be affiliated 
 with $A$, since it is well known that a sum and a product of two 
closed operators may not be closed. Affiliation relation in \cstal $\;$\  theory 
  was introduced by   Baaj and  Julg in  \cite{baajaf}. 

Observe, that if  $\phi\in \Mor (A,B)$, 
then one can extend  $\phi$ 
 to elements  affiliated with  $A$. Let us start with the  observation, 
 that for any 
$T\in M(A)$ we have 
\[ \phi (z_T)=z_{\phi (T)}\ .\]
 Hence for any  $T\af A$ we have  $z_T\in M(A)$. Moreover, there exists   
  a unique closed  operator S such that  $\phi (z_T)=z_S$. 
This operator is given by  
$$S=\phi (z_T)\phi (I-z_T^*z_T)^{-\frac{1}{2}}.$$
 From now on we will  write  
 $S=\phi (T)$. 

 We recall now a  nonstandard notion of  generation we use in this paper.  
This notion was introduced in  \cite{wunb}, where  a  
 generalization of the theory of unital  \cstal s generated by  
 a finite number of generators was presented. It was proved in 
\cite{oper} that such \cstal s  
 are isomorphic to  algebras of  
 continuous operator functions on compact  operator domains
(see  Section 1.3 of  \cite{paper1} and references therein). 
In this approach, the  algebra of all  continuous vanishing at infinity 
 functions on a compact quantum    group  is generated  by matrix 
elements of  fundamental representation. 
  To use this approach to non compact quantum  groups,
  one has to extend  the notion of  a generation of a  \cstal $\;$\  to   
 non unital \cstal s  and  unbounded generators . 
According to the definition we recall below,  
{\em \cstal $\;$\ of continuous  vanishing at infinity functions 
 on a locally  compact quantum  group}
 is  generated  by its fundamental representation . 
However, in this case the fundamental  representation is not 
  unitary and  the generators are unbounded  operators,  
 so they are not in the \cstal $\;$\  $A$.  

Assume for a while, that were are  given  a \cstal $\;$\  $A$  and operators 
 $T_1,T_2,...,T_N$ affiliated with   $A$.
We say that $A$ is {\em  generated} by  $T_1,T_2,...,T_N$ 
if for any Hilbert space 
$\cH$,  a non degenerate \cstal $\;$\ $B\subset B(\cH)$ 
and any $\pi\in Mor(A,CB(\cH))$ we have 
\[\left(
\begin{array}{c} 
\pi (T_i) \mbox{ is affiliated with } A\\
\mbox{ for any } i=1,...,N  
\end{array}
\right)
\Longrightarrow
\left(
\begin{array}{c}
\pi\in \Mor(A,B)
\end{array}
\right).\] 

We stress that described above 'generation' is a relation 
between $A$ and some operators $ T_1,T_2,...,T_N$ and 
both have to be known in advance. There is no 
procedure to obtain $A$ knowing only $ T_1,T_2,...,T_N$ and 
it is even possible that  there is no  $A$ generated by such operators.

For unital  \cstal s  generation in the sense introduced above 
 is exactly the same as the classical notion of generation. 
More precisely, 
  let  $A$ be a unital  \cstal $\;$   \ and let  
\mbox{$T_1, T_2, \dots , T_N\in A$.}
If  $A$ is the norm closure of all linear  combinations of  
\mbox{$ I, T_1,  \dots , T_N$,}
 then  $A$ is generated  by $ T_1, T_2, \dots , T_N$ in the sense of the  
 above  definition. 
On the other hand, let    $A$ be a  \cstal $\;$\ generated  by 
 $T_1, T_2, \dots , T_N\af A$, such that  $||T_i||<\infty $ for 
$i=1,2,\dots, N$. Then $A$ contains unity , 
$T_1, T_2, \dots , T_N\in A$ and  $A$ is the norm closure of the set of all 
  linear  combinations of 
$ I, T_1, T_2, \dots , T_N$.

An easy example of this relation is that \csta algebra $\Cir$ 
of all continuous vanishing at infinity functions on $\R$ 
is generated by function $f(x)=x$ for any $x\in \R$.
The other example is \cstal $\;$ \ $\zw$ which is generated by 
momentum operator $\po$ and multiplication-by-coordinate operator  $\qo$.

Let   $A$ and $B$ be  \cstal s and assume that we know generators of $A$. 
In order to describe  $\phi \in \Mor(A,B)$ uniquely it is enough 
  to know how  $\phi$ acts on  generators of  $A$.

We will use exclusively the minimal tensor product of \cstal s  and it will 
be denoted by $\te$.
We will also use the leg numbering notation. For example, if $\phi
\in \Mor (A\te A,A\te A)$ then $\phi_{12}(a\te b)= a\te b\te I_A$
and  $\phi_{13}(a\te b)= a\te I_A\te b$ for $a,b\in A$.
Clearly, $\phi_{12}, \phi_{13}\in \Mor (A\te A,A\te A\te A)$.

Let  $f$ and $\phi$ be strongly commuting selfadjoint operators.
Then, by the spectral theorem
\[f=\int _{\Lambda} r dE(r,\rho)\mi \phi=\int _{\Lambda} \rho
 dE(r,\rho)\ ,\]
where $dE(r,\rho)$ denotes the common spectral measure associated 
with $f$ and $\phi$ and $\Lambda$ stands for a joint spectrum 
of $f$ and $\phi$.
Then 
\[
F (f,\phi)=\int _{\Lambda} F(r,\rho) dE(r,\rho)\ .
\]

Let $b$ be a selfadjoint operator and  let  the symbol $\chi$ 
denote the characteristic function defined on  $\R$. 
By $\chi(b\neq 0)$ we mean  the projection  operator on 
 the subspace $\ker b^\perp$, by $\chi(b< 0)$ - the projection  onto  
 the subspace on which $b$ is negative, and so on.

\section{Relevant information on the quantum ``az+b'' group}
\label{qa}

Let
\beq
\label{q}
q=e^{\frac{2\pi i}{N}}\ ,
\eeq
where $N$ is an even number and  $N\geq 6$, 
i.e. $q$ is a primitive root of unity: $q^N=1$.
Let us introduce notation
\beq
\label{hb2}
\hb=\frac{2\pi}{N}\ .
\eeq
Note that $\hb<\pi$ and  $q=e^{i\hb}$.

Let
\beq
\label{gamma}
\Gamma =\bigcup_{k=0}^{N-1}q^k\;\R_+\ .
\eeq
Observe that  $\Gamma$ is a multiplicative subgroup of $\C\setminus 
\{0\}$. 

Let  $\ov\Gamma$ denote a closure of  $\Gamma$ in $\C$, i.e.
\beq
\label{gammaov}
\ov\Gamma =\Gamma \cup \{0\}\ .
\eeq

We introduce now  the  definition of the operator domain 
  $G_{\sH}$ \cite{paper3,az+b}, related to the quantum    "az+b" 
group (for the notion of operator domain we refer the Reader to 
\cite{paper1, phd} and references therein)

\[ G_{\sH}=\left\{(a,b)\in\ch^2\;|\; \ba{c}
   aa^*=a^* a,\\  bb^*=b^*b,\\                                                  
\ker a=\{0\}\\        
 \spe a, \spe b\subset\ov\Gamma\\
      (\faz a )b=qb (\faz a)\\
\mbox{\rm for  any  } t\in\R\\
|a|^{it}b=e^{-\frac{2\pi}{N}}b|a|^{it}
\ea
\right\}\]     

Observe that the relationship between $G_{\cH}$ 
and braided quantum group  $D_{\cH}$ considered in 
\cite{paper3,phd} is the following
\[\left( \ba{c} (b,a)\in D_{\cH}^2\\ \ker a=\{0\} \ea \right) 
\Longleftrightarrow (a,b)\in G_{\cH}\]

Let   $A$
denote the  \cstal $\;$    of all continuous functions vanishing at infinity on 
 the quantum  ``az+b'' group. The algebra $A$   
 is   generated  (in the sense explained in Section
  \ref{not})
  by unbounded operators  $a$, $\;$ $a^{-1}$ i $ b$, 
where  $(a,b)\in G_{\sH}$.
It was proved in  \cite{az+b} that 
the  multiplicative unitary operator 
\footnote{To prove  manageability of  $W$ we need a slightly 
 more complicated 
formula \cite{az+b}. However, for the purpose of this paper, we 
may consider a simpler formula  \rf{W}.}
 $W\in B(\cH\te \cH)$\index{$W\;\;$} 
 for the \qa is given by
\beq
W=F_N \left(a b^{-1} \te b\right)
\chi(b^{-1}\te I,I\te a)\label{W}
\ ,\eeq
where $\chi$ is a  symmetric  bicharacter on $\Gamma$, such that 
 for any  $\gamma\in\Gamma$ and $r\in\R_+$ 
\[ \chi (\gamma, q)=\faz \gamma \mi \chi (\gamma,r)=|\gamma|^{\frac{N}{2\pi i}\log r}.\]

From the  theory of   multiplicative  unitaries  we know that operator 
 $W$ encodes the group structure of the quantum group.
 Precisely, for any  $d\in A$ comultiplication 
 \mbox{ $\Delta\in \Mor (A,A\te A)$} 
is given by
\[\Delta(d)=W(d\te \id)W^*\ .\]
Comultiplication  $\Delta$ may be extended  to unbounded 
 operators affiliated with 
 $A$ and is given on generators of $A$  by the same formula
\[\Delta(a)=W(a\te \id)W^*\]
\[\label{deltab}\Delta(b)=W(b\te \id)W^*\]
It was also computed in \cite{az+b} that for any  $(a,b)\in G_{\sH}$ 
we have
\[\label{deltaa}\Delta(a)=a\te a\]
\[\Delta(b)=a\te b\dot{+}b\te I\ ,\]
where $a\te b\dot{+}b\te I$ denotes the closure of the sum   $a\te b+b\te I$.

Thus  defined $\Delta$ is associative.
Moreover, $G$ equipped with thus defined  $\Delta$ is the \qa.
However, what we are mostly interested in  in this paper  is 
the  \cstal $\;$  $A$   of all continuous functions vanishing at infinity on 
 the quantum  ``az+b'' group.

We remind now the  construction of a certain   \csta 
 dynamical system investigated in \cite{az+b}. It turns out that 
 the  \cstal $\;$ corresponding to this system is 
 exactly our algebra $A$.

Let 
$$B=\{f\in \Ci(\ov{\Gamma})\}$$
and let  
$$b(\gamma)=\gamma$$
for any  $\gamma\in \ov{\Gamma}$.
Then  $b\af B$ and  $b$ generates  $B$ in the sense of  \cite{wunb}.

Let us define an action 
$\sigma\in {\rm Aut}(\M (B))$\index{$\sigma\;\;$} of the group   
${\Gamma}$ on any  function  $f\in B$ 
by 
\[(\sigma_t f)(\tau)=f( t\tau)\whe \tau\in\ov{\Gamma}
 \mi t\in {\Gamma}\ . \]
Then $(B,{\Gamma},\sigma)$ is a   \csta   dynamical 
system  
 (see e.g. \cite{ped}). 
Denote the    \csta - crossed product algebra 
 by  $A_{cp}$
\[A_{cp}=B\times_{\sigma} {\Gamma} \ .\]
By definition of the \csta - crossed product algebra,  
 $M(A_{cp})$
contains 
a strictly continuous  one-parameter group of  unitary 
operators,  
 implementing an action  $\sigma$ of the group  ${\Gamma}$ on an  algebra $B$:
\[\sigma_t f=U_tf U_{-t}\ \]
for any $f\in\B$ and $\gamma\in \Gamma$.
By  SNAG Theorem and  Theorem  5.2 \cite{az+b}, every  
 strictly continuous  one-parameter group of  unitary operators 
$(U_t)_{t\in\Gamma}$ contained in  $M(A_{cp})$ has form 
\[U_t=\chi (a,t),\]
where  $ a$ is a normal operator affiliated  with $A_{cp}$
 and $\spe a\subset \ov{\Gamma}$. Moreover, $a$ is invertible and  
$a^{-1}\af A_{cp}$.
Moreover, it can be easily seen that 
\[\sigma_t b=t b\ .\]
Hence
\[U_t b=t b U_t,\]
for any  $t\in\Gamma$.
It follows that 
$(a,b)\in G_\cH$.

It is also well-known that the linear envelope of a set  
\[ \{\ fg( a)\ : f\in B, \ g\in\Ci (\Gamma)\}\]
 is dense in   $A_{cp}=B\times_\sigma{\Gamma}$. 
Hence $B\subset M(A)$.
It turns out (\cite[Proposition 4.1]{az+b}) that  
 $a,a^{-1} $ and  $b$ generate  $A_{cp}$ in the   sense of   \cite{wunb}.
Moreover,  Proposition 3.2 \cite{az+b} says that for any pair 
$(\tilde{a},\tilde{b})\in G_{\cH}$ there is unique  representation  
$\pi\in {\rm Rep} (A_{cp},H)$ such that  $\tilde{a}=\pi (a)$ and   
 $\tilde{b}=\pi(b)$. 
It was proven in \cite{az+b} that  $A_{cp}=A$. 
From know on we will not distinguish $A_{cp}$ and $A$.

We proceed to construct a dual action of the group  ${\Gamma}$ on   
  $A=B\times_{\sigma} {\Gamma} $.
To this end we consider a map
\[
\label{teta}
\theta_{\gamma} =(\id\te\fil _{\gamma})\Delta \ ,
\index{$\theta_{\gamma}\;\;$}\]
 where a map 
$\fil \in \Mor (A,\Ci (\Gamma))$ is such that  
 for any   $t\in{\Gamma}$
\[
\fil _ {\gamma}( a)=\gamma \mi  \fil _{\gamma} (b)=0\ .
\label{zgen}
\]
Proposition 4.2   \cite{az+b} implies that there is  only one such 
a map  $\phi$.
Observe that    $\theta_\gamma$ is an   automorphism of $A$ and that  
\[\theta_0=\id \mi \theta_{\gamma_1}\theta_{\gamma_2}=
\theta_{\gamma_1\gamma_2}\ \]
for any  $\gamma_1,\gamma_2\in\Gamma$.
The map $\gamma\rightarrow \theta_\gamma(d)$ is 
continuous for any 
  $d\in A$.

Moreover, for any  $\gamma,t\in {\Gamma}$
\[\theta_\gamma(U_{t})=\chi(\gamma,t) U_t\ . \]

Thus  we showed that  $(A,{\Gamma},\theta)$ is a 
\csta dynamical system and is dual  
 to the  \csta dynamical system $(B,{\Gamma},\sigma)$.

Using  \rf{deltaa} and \rf{zgen} we compute
\[
\label{tetaa}
\theta_\gamma(a)= (\id\te\fil _\gamma)\Delta (a)=
(\id\te\fil _\gamma)(a\te a)=\gamma a\ \]
and 
 \[
\label{tetab}
\theta _\gamma (b)=(\id\te\fil _\gamma)\Delta (b)=b\ .\]
Hence for any   $\gamma\in \Gamma$ and 
$g\in \M(B)$ we have 
\[\theta_\gamma (g)=g \ .\]
Since $(B,{\Gamma},\sigma)$ and  $(A,{\Gamma},\theta)$ are  dual 
 \csta-dynamical systems, it follows that  
$(B^{\prime\prime},{\Gamma},\sigma)$ and 
$(A^{\prime\prime},{\Gamma},\theta)$ are  dual $W^*$-dynamical systems,
 if only  $\sigma$ and  $\theta$ are extended in the obvious way. 

Let  $\cK$ be a Hilbert space  (this time we consider  
also finite-dimensional ones).
Moreover, observe that it follows from the above remark that 
 also 
$(B(\cK)\te B^{\prime\prime},{\Gamma},I\te\sigma)$
is a  $W^*$-dynamical system  and its  
  von Neumann $W^*$- crossed product algebra is  
$B(\cK)\te A^{\prime\prime}$.
Analogously, the   $W^*$-dynamical system dual to   
  $(B(\cK)\te B^{\prime\prime},{\Gamma},I\te\sigma)$ 
 is   $(B(\cK)\te A^{\prime\prime},{\Gamma},I\te\theta)$.

In what follows we need Proposition \ref{niezm}, 
which is an obvious consequence 
of  Theorem  7.10.4   \cite{ped}:
\bst
\label{niezm}
Let  $m\in B(\cK)\te A^{\prime\prime}$ and let for any 
 $\gamma\in {\Gamma}$  
\[ (\id\te \theta_\gamma)(m)=m\ .\]
Then
\[m\in B(\cK)\te B^{\prime\prime}\ .\]
\est
We also need
\bst
Let   $B$ be a     \cstal $\;$ and  let  $w\in\M(B\te A)$. 
Then a map   $\ov{\Gamma}\ni \gamma
\rightarrow (\id\te\fil_\gamma) w \in \M(B)$ 
 is strictly continuous.
\label{dostone}
\est
\dow We know that    $\fil \in \Mor(A,\Ci ({\Gamma}))$
 and therefore 
  $(\id \te \fil)w\in 
\M(B\te \Ci ({\Gamma}))$. S.L. Woronowicz  proved in  \cite{wunb} 
that 
  elements of 
$\M(B\te \Ci ({\Gamma}))$ are bounded, strictly 
   continuous   functions    
on  ${\Gamma}$ with values in   $M(B)$.
\hfill\qed

\section{How do  unitary representations of the quantum ``az+b'' 
group look like?}

\label{ost}

\bde
A unitary operator  $V\in \M(CB(\cK)\te A)$ is called a    
 (strongly continuous) {\em unitary representation}  of the \qa 
 if 
\[W_{23}V_{12}=V_{12}V_{13}W_{23}\ ,\] 
or equivalently
\beq
\label{ru}
(\id\te \Delta)V=V_{12}V_{13}\ .
\eeq
\ede

Observe that in case of the classical group  condition  \rf{ru}
 is equivalent to the usual  definition of a  unitary representation, i.e.  
a representation  $U$  is a map 
$$U:G\ni g\rightarrow U_g\in B(\cK)$$
such that  $U_g$ is  unitary  for any  $g\in G$ and 
 for any $g,h\in G$ we have $U_g U_h=U_{gh}$. 

It was proven in  \cite{az+b} that
\bst
Let   $(a,b)\in G_{\sH}$ and  $(c,d)\in G_{\sK}$ and let 
$\ker b=\{0\}$.
Then the operator $V\in \M(CB(\cH)\te A)$ given by 
\[
\label{V}
V(a,b)=F_N(d \te b)
\chi(c \te I,I\te a)
\] 
is a  unitary representation  of the \qa
\est
We now prove  that all  unitary  representations  of the \qa 
are of this form.
.
\btw
Every   unitary  representation of the \qa (i.e. of  $G$), 
acting on  some 
 Hilbert space   $\cK$ 
 has form 
\[V(a,b)=F_N(d \te b)
\chi(c \te I,I\te a)
\] 
where   $(c,d)\in G_{\sK}$.
\etw
\dow
Let   $V$ be a  unitary representation of the \qa
acting on  a 
 Hilbert space   $\cK$. Then for any 
$\gamma\in {\Gamma} $ 
 an operator $(\id\te\fil _\gamma)V\in B(\K)$ is unitary.
Applying    $(\id\te\fil_s\te\fil_t)$ to both sides of   
 \rf{ru} we get 
\[V(\hb (s+t),0,0)=V(\hb s,0,0)V(\hb t,0,0)\]
Hence
\[(\id\te\fil_{s+t})V=(\id\te\fil_{s})V(\id\te\fil_{t})V\]
i.e.   $(\id\te\fil)V$ is a   representation of ${\Gamma}$.

The strict topology on  $B(H)=M(CB(H))$ coincides with 
 the  *-strong   operator topology.
Since  the  *-strong   operator topology is stronger than 
 the strong   operator topology, then by  
  Proposition   \ref{dostone}  we obtain that the map 
\[ {\Gamma}\ni t\rightarrow (\id\te\fil_t)V\in B(\K)\]
 is strongly continuous.
 Therefore by the  SNAG Theorem there is a 
 normal operator  $c$ acting on    $\K$ 
 with spectrum contained in $\Gamma$ and such that 
\[
\label{zau1}
(\id\te\fil_\gamma)V=\chi (c,\gamma)\]
 for any   $t\in{\Gamma}$.

Note that from \rf{teta}, \rf{ru} and  \rf{zau1} follows that 
\[(\id\te\theta_\gamma)V=V\chi( c\te I,\gamma I\te I))\ .\]
Moreover, by 
\[(\id\te\theta_\gamma)\chi(c\te I, I\te a)^*
=\chi(c\te I,I\te \gamma a)^*=
\chi(c\te I,\gamma I\te I )^*\chi(c\te I,I\te  a)^*
\ .\]
Hence 
\[(\id\te\theta_\gamma)V\chi(c\te I, I\te a)^*=
V\chi(c\te I, I\te a)^*\ .\]
Observe that 
\[V\chi(c\te I, I\te a)^*\in B(\cK)\te A^{\prime\prime}\ .\]
Hence by Proposition  \ref{niezm}
\[V\chi(c\te I, I\te a)^*=f(b)\ ,\]
where  
$f\in  B(\cK)\te B^{\prime\prime}$, i.e.  $f$ is a Borel operator 
 function  on  $ \Gamma$ with values in bounded operators acting on    \
$\cK$ (for explanation on operator functions see 
\cite{paper1}[Section 1.3] and references therein). 
Hence
\[V=f(b)\chi(c\te I, I\te a)\ .\]
Since     $V$   is unitary   operator  , 
it follows from the above formula that 
 \mbox{ $f(b)\in {\rm Unit}(\K\te\cH)$.} 
Compute 
\[(\id\te\Delta)V=(\id\te\Delta)\left(f(b)
\chi(c\te I, I\te a)\right)=f(\Delta(b))
\chi(c\te \Delta(I), I\te \Delta(a))=\]
\[=f(a\te b\dot{+}b\te I)\chi(c\te I\te I, I\te a\te a)\ .\]
Applying    $(\id \te \fil_\gamma\te \id)$ to both sides of \rf{ru}
 we obtain, respectively,
\[(\id \te \fil_\gamma\te \id)V_{12}V_{13}=(\chi(c,\gamma)\te \id)V\]
and
\[(\id \te \fil_\gamma\te \id)(\id\te\Delta)V=f(\gamma b)
\chi (c\te I,\gamma\te a)\ .\]
Comparing these results we see that 
\beq
\label{por1}
(\chi(c,\gamma I)\te \id)V
=f(\gamma b)
\chi (c\te I,\gamma I\te a)\ .
\eeq
On the other hand, applying    
$(\id \te\id \te \fil_\gamma)$ to both sides of \rf{ru}
 we conclude that 
\[(\id \te \id \te \fil_\gamma)V_{12}V_{13}=V(\chi(c,\gamma I)\te \id)\]
and 
\[(\id \te \id \te \fil_\gamma)(\id\te\Delta)V=f(b)
\chi (c\te I,\gamma I\te a)
\ .\]
Comparing these results we see that 
\beq
\label{por2}
V(\chi(c,\gamma I)\te \id)=f(b)
\chi (c\te I,\gamma I\te a)\ .
\eeq
Inserting             $ a$ in the place of $\gamma$ 
 in formulas  \rf{por1} and \rf{por2} we derive 
\[(\chi (c\te I,I\te a)\te \id)V_{13}
=f(a\te b)\chi (c\te I\te I,I\te a\te a)
\]
and
\[V_{12}(\chi (c\te I,I\te a))_{13}
=f(b\te I)
\chi (c\te I\te I,I\te a\te a)\ .\]
Hence
\[V_{12}=f(b\te I)\chi (c\te I\te I,I\te a\te I)\]
\[V_{13}= \chi (c\te I\te I, I\te a\te I)^* f(a\te b)
\chi (c\te I\te I, I\te a\te a)
\label{565}\ .\]
Therefore, since   $V$ satisfies 
  $(\id \te\Delta)V= V_{12}V_{13}$, we have 
\[
f(a\te b\dot{+}b\te I)\chi ( c\te I\te I, I\te  a\te a)=\]
\[=f(b\te I) f(a\te b)
\chi ( c\te I\te I, I\te  a\te a)
\ .\]
Hence obviously
\beq
\label{dlatego}
f(a\te b\dot{+}b\te I)
=f(b\te I) f(a\te b)
\ .
\eeq
Let us introduce notation  
\beq
R=a \te b\  \od \label{oznrs}
\mi
S= b\te I\ \ .
\eeq
We conclude that    $(S,R)\in D^2$. 
Inserting notation \rf{oznrs}
 into formula  
  \rf{dlatego} we derive 
\[
f(S \dot{+}R)
=f(S)f(R).\]
By Theorem 7.1 \cite{paper3} if a    function  $f$ is a Borel 
operator function  $f:D\rightarrow B(\cK\te\cH)$ and satisfies 
 the above condition for $(S,R)\in D^2$, then it is given by
\[ f(b)=F_N \left( d\te b\right) \ ,\]
where $d$ is a normal operator  with spectrum contained in  $\ov\Gamma$.

What is left is to prove that  $(c,d) \in G_{\cK}$.
To this end observe that by   \rf{565} and from 
\[V=F_N(d\te b)\chi (c\te I,I\te a),\]
follows that 
\[\fh (d\te I\te b)\chi
(c\te I\te I,I\te I\te a)= \]
\[=\chi
(c\te I\te I,I\te a\te I)^*
F_N(d\te a\te b)\chi
(c\te I\te I, I\te a\te a)\ .\]
Hence
\[ \label{stad1} (d\te \id)=
\chi (c\te I,I\te a)^*
(d\te a) 
\chi (c\te I, I\te a)\ .\]
Inserting    $a=\gamma I$ in the formula  \rf{stad1} we obtain that  
   $(c,d,\delta)\in G\sK$, which completes the proof.\hfill\qed


\begin{thebibliography}{99}

\bibitem{baajaf} S. Baaj \& P. Julg, "Th\'{e}orie bivariant de Kasparov et 
op\'{e}rateur non born\'{e}s dans les \csta modules hilbertiens", 
{\em C. R. Acad. Sci. Paris}, S\'{e}rie I, {\bf 296} (1983) 875-876.  

\bibitem{ks}
Klimyk, Schm\"{u}dgen: Quantum Groups and Their Representations.Springer
-Verlag Berlin -- Heidelberg 1997.

\bibitem{oper}
P. Kruszy\'nski \& S.L. Woronowicz: A Non-commutative Gelfand-Naimark Theorem-
J. Operator Theory 8 (1982), 361-389. 


\bibitem{vaesk}  J. Kustermans, S.Vaes, A simple definition for locally compact quantum group {\em C.R. Acad. Sci. Paris, 
 S\'er. I } {\bf 328} (10)(1999), 871-876.


\bibitem{ped}G. K. Pedersen, {\em \cstal s and Their Automorphism Groups},
 Academic Press, London, New York, San Francisco, 1979. 


\bibitem{phd}
M. Rowicka: PhD Thesis, Warsaw University, Warsaw 2000.  

\bibitem{paper1}
M. Rowicka: Braided quantum groups related to the quantum ``ax+b'' group 
  -  math.QA/0101003.

\bibitem{paper2}
M. Rowicka: Unitary representations of the quantum 'ax+b' group 
- math.QA/0102151.

\bibitem{paper3}
M. Rowicka: Exponential equations for the   quantum 'az+b' group 
-  math.QA/0103086.



\bibitem{nowywandal} van Daele, A.: The Haar measure on some locally compact 
quantum groups -in preparation.

\bibitem{qexp}
S.L. Woronowicz: Quantum exponential function -Rev. Math. Phys.
 Vol. 12, No. 6 (2000) 873-920.

\bibitem{ax+b}
S.L. Woronowicz \& S. Zakrzewski: Quantum 'ax+b' group - submitted to 
Comm. Math. Phys.

\bibitem{oe2}
S.L. Woronowicz: Operator Equalities Related to the Quantum E (2) Group - 
Commun. Math. Phys. 144, 417-428 (1992).

\bibitem{wunb}
S.L. Woronowicz: \csta algebras generated by unbounded elements.
{\em Rev. Math. Phys.} Vol.{7} No. 3 (1995) 481-521.
 
\bibitem{az+b}
Woronowicz, S.L.: Quantum 'az+b' group on complex plane -  KMMF Preprint 1999.

\end{thebibliography}
\end{document}